\documentclass[10pt]{article}
\usepackage{times}
\usepackage{amsmath,amsfonts,amssymb,latexsym,epsfig}
\usepackage{color,epsf}
\usepackage{graphicx}
\newtheorem{proposition}{Proposition}
\begin{document}
\title{Numerical integrators for the Hybrid Monte Carlo method}
\author{
Sergio Blanes\footnote{Instituto de Matem\'atica Multidisciplinar,
  Universitat Polit\`ecnica de Val\`{e}ncia, E-46022  Valencia, Spain.
 {\tt serblaza@imm.upv.es}},
 Fernando Casas\footnote{Institut de Matem\`{a}tiques i Aplicacions de Castell\'o and
   Departament de Matem\`{a}tiques, Universitat Jaume I,
  E-12071 Castell\'on, Spain.
{\tt Fernando.Casas@uji.es}},
J.M.\ Sanz-Serna\footnote{Departamento de Matem\'atica Aplicada e IMUVA, Facultad de Ciencias,
  Universidad de Valladolid, Valladolid, Spain.
{\tt sanzsern@mac.uva.es}}
}
\maketitle

\begin{abstract}We construct numerical integrators for Hamiltonian problems that may advantageously replace the standard Verlet time-stepper within Hybrid Monte Carlo and related simulations. Past attempts have often aimed at boosting the order of accuracy of the integrator and/or reducing the size of its error constants; order and error constant are relevant concepts in the limit of vanishing step-length.
We propose an alternative methodology based on the performance of the integrator when sampling from Gaussian distributions with not necessarily small step-lengths. We construct new splitting formulae that require two, three or four force evaluations per time-step. Limited, proof-of-concept numerical experiments suggest that the new integrators may  provide an improvement on the efficiency of the standard Verlet method, especially in problems with high dimensionality.
\end{abstract}\bigskip

\noindent AMS numbers: 65L05, 65C05, 37J05

\noindent Keywords: Hybrid Monte Carlo method, Markov Chain Monte Carlo, acceptance probability, Hamiltonian dynamics, reversibility, volume preservation, symplectic integrators, Verlet method, split-step integrator, stability, error constant, molecular dynamics

\section{Introduction}

The present paper constructs numerical integrators for Hamiltonian problems that may advantageously replace the standard Verlet time-stepper within Hybrid Monte Carlo (HMC) and related simulations.
HMC, introduced in the physics literature by Duane et al. \cite{duane}, is a Markov Chain Monte Carlo method \cite{robert} that has the potential of combining global moves with high acceptance rates, thus improving on alternative techniques that use random walk proposals \cite{neal}, \cite{natesh}. It is widely used in several areas, including quantum chromodynamics \cite{joo}, \cite{takaishi}, and is becoming increasingly popular in the statistics literature as a tool for Bayesian inference (see e.g. \cite{neal}). A longer list of references to various application areas may be seen in \cite{natesh}. At each step of the Markov chain, HMC requires the numerical integration of a Hamiltonian system of differential equations; typically, the familiar Verlet algorithm \cite{schlick} has been used to carry out such an integration. Since the bulk of the computational effort in HMC lies in the simulation of the Hamiltonian dynamics, it is of clear interest to investigate whether the simple Verlet algorithm may be replaced by more sophisticated and efficient alternatives. In particular better integrators may reduce the number of rejections, something valuable in applications such as molecular dynamics where discarding a computed trajectory may be seen as a significant \lq waste\rq\ of computational time. Although the physics literature is not lacking in efforts to construct new integrators (see e.g.\ \cite{joo}, \cite{takaishi} and their references), the fact is that Verlet remains the integrator of choice.

Past attempts to build integrators to improve on Verlet have typically started from the consideration of families of split-step methods with one or several free parameters; the values of those parameters are then adjusted to boost the order of accuracy  and/or to reduce the size of the error constants. We shall argue here that such a methodology, while well-established in numerical analysis, cannot be expected to be fruitful within the HMC context. In fact, order of accuracy and error constants are notions that provide information on the behavior of an integrator as the step-size $h$ approaches 0 and  in HMC simulations useful integrators operate with moderate or even large values of $h$. In an alternative approach, we begin by associating with each numerical integrator a quantity $\rho(h)$ that governs its behavior  in  simulations of Gaussian distributions (Proposition \ref{prop:average} and Section \ref{ss:multivariate}). More precisely $\rho(h)$ provides an upper bound for the energy error when integrating the standard harmonic oscillator and is relevant to all multivariate Gaussian targets. We then choose the values of the free parameters to minimize the size of $\rho(h)$ as $h$ ranges over an interval $0<h<\bar{h}$, where $\bar{h}$ is sufficiently large.\footnote{This approach is somewhat reminiscent of the techniques used in \cite{01} and \cite{02}.} Numerical experiments show that the new approach does produce integrators that provide substantial improvements on the Verlet scheme. On the other hand, when  integrators derived by optimizing error constants and Verlet are used with step-lengths that equalize work, the energy errors of the former typically improve on those of Verlet only for step-sizes so unrealistically small that the acceptance rate for Verlet is (very close to) 100\%.

After submitting the first version of the present work, we have become aware of two additional references, \cite{nuevo} and \cite{benmath}, that are relevant to the issues discussed here. The paper \cite{nuevo} considers integrators for Hamiltonian dynamics and, just as in the present work, tunes the coefficients of the methods so as to ensure good conservation of energy properties in linear problems; furthermore \cite{nuevo} discusses the reasons the relevance of linear models as guides to nonlinear situations. However the optimization criterion of \cite{nuevo} differs from ours, as it is based on maximizing the length of the stability interval, subject to the annihilation of some error constants. In \cite{benmath} the authors deal with Langevin integrators and demonstrate methods which have exact sampling for Gaussian distributions.

Sections 2  and 5 provide the necessary background on HMC and  splitting integrators respectively; in order to cater for readers with different backgrounds the exposition there is rather leisurely. Section 3 studies a number of peculiarities of the numerical integration of Hamiltonian systems {\em specific to the HMC scenario}. We point out  that the average size of the energy error is actually much smaller than one would first believe.  In such a scenario the optimal stability property of Verlet makes the construction of a  more efficient integrator a rather demanding challenge. Our methodology for determining the free parameters in families of integrators is based on Gaussian model problems; such models are studied in Section 4. We show in particular that, for Gaussian targets and if the dimensionality is not extremely large, the Verlet algorithm performs well with values of the step-length $h$ that are moderate or  large. Section 6 presents our approach to the choice of free parameters. It also contains examples of methods with two, three or four force evaluations per time-step derived by following the new methodology. The new methods clearly outperform the Verlet integrator, particularly so if the dimensionality of the problem is high. Section 7 reports some numerical comparisons in a simple molecular example and the final Section 8 is devoted to conclusions.

\section{The hybrid Monte Carlo method}

The aim of the HMC algorithm is to obtain a Markov chain \cite{robert} to sample from a probability distribution in $\mathbb{R}^d$ with density function of the form
\begin{equation}\label{eq:pdf}
\pi(q) \propto \exp(-V(q)).
\end{equation}
The algorithm introduces an auxiliary variable $p\in\mathbb{R}^d$, called {\em momentum}, and works in the {\em phase space} $\mathbb{R}^{2d}$ of the variables $(q,p)$, where
one considers a {\em Hamiltonian} function $H$ (energy)
\begin{equation}\label{eq:hamiltonian}
H(q,p) = \frac{1}{2} p^TM^{-1}p+V(q)
\end{equation}
($M$ is a symmetric, positive definite matrix chosen by the user)
and a probability density function
\begin{equation}\label{eq:pdfpq}
\Pi(q,p) \propto \exp\big(-H(q,p)\big) = \exp\big(-\frac{1}{2} p^TM^{-1}p\big)\:\exp\big(-V(q)\big).
\end{equation}
Thus $q$ and $p$ are stochastically independent,  $q$ is distributed according to the target (\ref{eq:pdf})
and $p$ has a Gaussian $N(0,M)$ distribution.
\begin{figure}
\begin{center}
\includegraphics[scale=0.35]{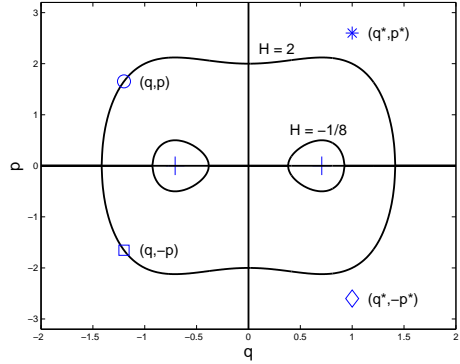}
\end{center}
\caption{\em The level sets $H=2$ and $H=-1/8$ in the phase plane $(q,p)$ when $V(q) = q^4-q^2$. The target distribution $\pi(q)$ has modes at $q = \pm \sqrt{2}/2$. A reversible transformation $\Psi$ that maps $(q,p)$ into $(q^*,p^*)$ must map $(q^*,-p^*)$ into $(q,-p)$. Here the move from $(q,p)$ into $(q^*,p^*)$ increases the value of $H$ and the move from
$(q^*,-p^*)$ to $(q,-p)$ decreases $H$ by the same amount.}
\label{fig:rever}
\end{figure}
\begin{table}
\hrule
\bigskip

Given $q^{(0)}\in\mathbb{R}^d$, $N \geq 1$, set $n=0$.

\begin{enumerate}
\item Draw $p^{(n)}\sim N(0,M)$. Compute $(q^*,p^*) =\Psi(q^{(n)},p^{(n)})$ ($q^*$ is the proposal).

\item Calculate $a^{(n)} = \min\big(1, \exp(H(q^{(n)},p^{(n)})-H(q^*,p^*))\big)$.

\item Draw $u^{(n)} \sim U(0,1)$. If $a^{(n)}>u^{(n)}$, set $q^{(n+1)}=q^*$ (acceptance); otherwise set $q^{(n+1)} = q^{(n)}$ (rejection).

\item Set $n= n+1$. If $n = N$ stop; otherwise go to step 1.
\end{enumerate}
\hrule
\caption{\em Basic HMC algorithm.  $M$ and $\Psi$ are respectively a positive definite matrix and a volume-preserving and reversible transformation in phase space; both are chosen by the user. The function $H$ is given in (\ref{eq:hamiltonian}).
 The algorithm generates a Markov chain $q^{(0)}\mapsto q^{(1)} \mapsto \dots \mapsto q^{(N)}$ reversible with respect to the target probability distribution (\ref{eq:pdf}).}
 \label{tab:alg1}
\end{table}

The algorithm uses transitions in phase space  $(q^{(n)},p^{(n)}) \rightarrow (q^{*},p^{*}) =\Psi(q^{(n)},p^{(n)})$  obtained through a mapping
$\Psi:\mathbb{R}^{2d}\rightarrow \mathbb{R}^{2d}$ that is {\em volume preserving,}
\begin{equation}\label{eq:volpre}
\det(\Psi^\prime(q,p)) = 1
\end{equation}
($\Psi^\prime$ is the Jacobian matrix of $\Psi$), and {\em reversible} (see Fig.~\ref{fig:rever}),
\begin{equation}\label{eq:reversible}
\Psi(q,p) = (q^*,p^*) \Leftrightarrow \Psi(q^*,-p^*) = (q,-p).
\end{equation}
If $S$ denotes the mapping in phase space $S(q,p) = (q,-p)$ (momentum flip), then the reversibility requirement reads
$$
\Psi\circ S\circ \Psi = S
$$
($\circ$ denotes composition of mappings, i.e. $(S\circ \Psi) (q,p) = S\big(\Psi(q,p)\big)$), or, since $S^{-1} = S$,
\begin{equation}\label{eq:reversible2}
\Psi^{-1} = S\circ\Psi\circ S.
\end{equation}

Table~\ref{tab:alg1} describes the basic HMC procedure. A proof of the fact that the algorithm generates a Markov chain reversible with respect to the target probability distribution (\ref{eq:pdf}) may be seen in \cite{franceses} or \cite{cetraro} (see also Section 2.2 in \cite{cances}, which contains additional references). HMC is of potential interest because by choosing $\Psi$ appropriately it is possible to have
a proposal $q^*$ far away from the current location $q^{(n)}$ while at the same time having a large probability $a^{(n)}$ of the proposal being accepted; that is not the case for the random walk proposal in the standard Metropolis algorithm. (In Fig.~\ref{fig:rever}, if the circle and the star correspond to   $(q^{(n)},p^{(n)})$ and $(q^*,p^*)$ respectively, then
the current location is in the neighborhood of the mode at $q=-\sqrt{2}/2$ and the proposal is close to the other mode.)

It is in order to recall that the reversibility of a Markov chain with respect to a target distribution is not by itself sufficient to ensure the ergodic behavior that is required for the chain to yield trajectories that may be successfully used to compute averages: additional properties like irreducibility are necessary. The discussion of these issues is outside the scope of the present work and the interested reader is referred e.g.\ to \cite{cances} or \cite{franceses}.

Many variants and extension of the procedure in Table~\ref{tab:alg1} have been suggested in the literature, see, among others, \cite{elena}, \cite{beskos}, \cite{horo}, \cite{iza}, \cite{neal}. It is not our purpose here to compare the merit of the different variants of HMC or to compare HMC with other sampling techniques.

\subsection{Using Hamiltonian dynamics}
 A potentially interesting choice of transformation $(q^*,p^*) =\Psi(q^{(n)},p^{(n)})$ would be obtained by fixing a number $T>0$ and setting $q^*=q(T)$, $p^*=p(T)$, where $(q(t),p(t))$ is the solution of the system of differential equations
\begin{equation}\label{eq:hamsys}
\frac{d}{dt} q = \nabla_p H(q,p),\qquad \frac{d}{dt} p = -\nabla_q H(q,p)
\end{equation}
with initial condition $q(0) = q^{(n)}$, $p(0) = p^{(n)}$. In more technical words, in this choice, $\Psi$  coincides with the $T$-flow $\varphi_T^H$ of the Hamiltonian system (\ref{eq:hamsys}) \cite{arnold}, \cite{hlw}, \cite{benbook}, \cite{ssc}. By selecting $T$ suitably large, one then obtains a point $(q^*,p^*) = \varphi_T^H(q^{(n)},p^{(n)})$ away from $(q^{(n)},p^{(n)})$. Furthermore,
(\ref{eq:hamsys}) implies
$$
\frac{d}{dt}H(q(t),p(t)) = 0
$$
(conservation of energy), so that in step 2 of the algorithm
$H(q^{(n)},p^{(n)})-H(q^*,p^*) = 0$ and, accordingly, the probability of acceptance is $a^{(n)} =1$.

It is important to note here that, for each choice of $T$, $\Psi=\varphi_T^H$ satisfies the requirements in Table \ref{tab:alg1}. In fact, the preservation of volume as in (\ref{eq:volpre}) is a well-known result of the Hamiltonian formalism, see e.g.\ \cite{arnold}, \cite{hlw}, \cite{ssc}. Moreover the transformation $\Psi=\varphi_T^H$ is  reversible as in (\ref{eq:reversible}); this is checked by observing that $(\bar q(t),\bar p(t)) =(q(T-t),-p(T-t))$ is the solution of (\ref{eq:hamsys}) with initial condition $\bar q(0) = q^*$, $\bar p(0) = -p^*$ and that $(\bar q(T), \bar p(T)) = (q^{(n)},-p^{(n)})$.

Unfortunately the choice $\Psi = \varphi_T^H$ is unfeasible: in cases of practical interest (\ref{eq:hamsys}) cannot be integrated in closed form and it is not possible to compute
 $\varphi_T^H(q^{(n)},p^{(n)})$. HMC then resorts to transformations
$\Psi$ that {\em approximate} the true flow $\varphi_T^H$; more precisely $(q^*,p^*)$ is obtained by integrating (\ref{eq:hamsys}) with initial condition $(q^{(n)},p^{(n)})$ with a suitable numerical method. Not all methods can be considered: $\Psi$ has to satisfy the requirements (\ref{eq:volpre}) and (\ref{eq:reversible}). The well-known Verlet method, that we describe next, is at present the method of choice.

\subsection{The Verlet integrator}
If $h>0$ and $I\geq 1$ denote respectively the step-size and the number of time-steps, a {\em velocity Verlet} integration starting from $(q_0,p_0)$ may be represented as
\begin{equation}\label{eq:Psi}
\Psi(q_0,p_0) = (q_I,p_I)
\end{equation}
where $(q_I,p_I)$ is the result of the time-stepping iteration:
\begin{eqnarray}
p_{i+1/2} &=&p_i -\frac{h}{2} \nabla_qV(q_i),\nonumber \\
q_{i+1} &=& q_i +h M^{-1}p_{i+1/2},\nonumber \\
p_{i+1}  &=&p_{i+1/2} -\frac{h}{2} \nabla_qV(q_{i+1}), \qquad i = 0, 1,\dots, I-1.\label{eq:verlet}
\end{eqnarray}

For our purposes, the velocity Verlet algorithm is best seen as a {\em splitting} algorithm (see \cite{acta} and \cite{03}), where the Hamiltonian (\ref{eq:hamiltonian}) (total energy) is written as a sum $H=A+B$ of two partial Hamiltonian functions,
\begin{equation}\label{eq:kpot}A = (1/2)p^TM^{-1}p,\qquad B = V(q)
\end{equation}
that correspond to the kinetic and potential energies respectively.
 The Hamiltonian systems corresponding to the Hamiltonian functions $A$, $B$ are given respectively by
(cf.\ (\ref{eq:hamsys}))
$$
\frac{d}{dt} q = \nabla_p A(q,p) = M^{-1}p,\qquad \frac{d}{dt} p = -\nabla_q A(q,p) = 0
$$
and
$$\frac{d}{dt} q = \nabla_p B(q,p) = 0,\qquad
\frac{d}{dt} p = -\nabla_q B(q,p) = -\nabla_q V(q),
$$
and may be integrated in closed form. Their solution flows are respectively given by
$$
(q(t),p(t))= \varphi_t^A(q(0), p(0)),\qquad q(t) = q(0)+tM^{-1}p(0),\quad   p(t) = p(0),
$$
and
$$
(q(t),p(t)) = \varphi_t^B(q(0), p(0)),\qquad q(t) = q(0),\qquad p(t) = p(0)-t\nabla_q V(q(0)).
$$
Thus the Verlet time-step\footnote{In the numerical analysis literature it is customary to write \lq step\rq\ rather than \lq time-step\rq. Here
we  use \lq step\rq\ to refer to the Markov chain transitions and \lq time-step\rq\ to refer to the integration. Steps are indexed by the superindex  $n$ and time-steps by
the  subindex $i$.} $(q_i,p_i)\rightarrow (q_{i+1},p_{i+1})$ in (\ref{eq:verlet}) corresponds to a transformation
in phase space $(q_{i+1},p_{i+1}) = \psi_{h}(q_i,p_i)$ with
\begin{equation}\label{eq:comp1}
\psi_{h} = \varphi_{h/2}^B \circ \varphi_{h}^A\circ \varphi_{h/2}^B,
\end{equation}
and the transformation $\Psi=\Psi_{h, I}$ over $I$ time-steps (see (\ref{eq:Psi})) to be used in the algorithm in Table~\ref{tab:alg1} is given by the composition
\begin{equation}\label{eq:I}
\Psi =\Psi_{h, I}=\: \stackrel{I  \ {\rm  times}}{\overbrace{\psi_{h}\circ \psi_{h}\circ \cdots \circ \psi_{h}}}.
\end{equation}
Here $\psi_{h}$ is volume-preserving as a composition of volume-preserving Hamiltonian flows. Furthermore $\psi_h$
is reversible because $\varphi_{h/2}^B$ and $\varphi_{h}^A$ are both reversible {\em and} the right-hand side of (\ref{eq:comp1})
is a palindrome (see (\ref{eq:reversible2})):
\begin{eqnarray*}
\psi_{h}^{-1} &= & \big(\varphi_{h/2}^B\big)^{-1} \circ \big(\varphi_{h}^A)^{-1}\circ \big(\varphi_{h/2}^B)^{-1}\\
&=&\big(S\circ \varphi_{h/2}^B \circ S\big) \circ \big(S \circ \varphi_{h}^A \circ S\big)\circ \big(S\circ\varphi_{h/2}^B \circ S\big)\\
&=& S\circ \psi_{h} \circ S.
\end{eqnarray*}
It then follows that $\Psi_{h, I}$ is volume-preserving and reversible.
 Note that $\Psi_{h, I}$ is an approximation to the true solution flow at time $Ih$: $\Psi_{h, I}\approx \varphi_{Ih}^H$. Since
 $\varphi_{Ih}^H$ preserves energy exactly, the transformation $\Psi_{h, I}$ may be expected to preserve energy approximately, so that
 in Table~\ref{tab:alg1}, $H(q^{(n)},p^{(n)})-H(q^*,p^*) \approx 0$ leading to large acceptance probabilities.

 \begin{table}
\hrule
\bigskip

Given $q^{(n)}, p^{(n)}\in\mathbb{R}^d$, $I \geq 1$, $h>0$.

\begin{enumerate}
\item Set $q=q^{(n)}$, $p=p^{(n)}$, $i=0$.
\item Evaluate $\nabla_qV(q)$ and set $p = p-(h/2)\nabla_q V(q)$.

\item Set $q = q+h M^{-1}p$, $i= i+1$.

\item If $i < I$, evaluate $\nabla_qV(q)$ and set $p = p-h \nabla_q V(q)$, go to step 3. Otherwise go to step 5.

\item Set $q^*=q$, evaluate $\nabla_qV(q)$, set $p^* = p-(h/2)\nabla_q V(q)$ and stop.

\end{enumerate}
\hrule
\caption{\em Velocity Verlet algorithm to find $(q^*,p^*) =\Psi(q^{(n)},p^{(n)})$ in Table~\ref{tab:alg1}.}
 \label{tab:verlet}
\end{table}

Alternatively, the roles of $q$ and $p$ and those of the potential and kinetic energies may be replaced, to obtain the {\em position Verlet}
time-stepping \cite{schlick} (cf.~(\ref{eq:verlet})):
\begin{eqnarray*}
q_{i+1/2} &=&q_i +\frac{h}{2} M^{-1}p_i, \\
p_{i+1} &=& p_i -h \nabla_qV(q_{i+1/2}), \\
q_{i+1}  &=&q_{i+1/2} +\frac{h}{2} M^{-1}p_{i+1}, \qquad i = 0, 1,\dots, I-1.
\end{eqnarray*}
This is obviously a splitting integrator:
\begin{equation}\label{eq:comp2}
\psi_{h} = \varphi_{h/2}^A \circ \varphi_{h}^B\circ \varphi_{h/2}^A.
\end{equation}

The bulk of the work required to implement the Verlet velocity or position algorithms comes from the evaluation of the gradient $\nabla_q V$. In this connection it should be noted that the value $\nabla_q V(q_{i+1})$  in (\ref{eq:verlet}) coincides with the value to be used at the beginning of the subsequent $i+1\rightarrow i+2$ time-step. Thus both the velocity and position versions  require essentially one evaluation of $\nabla_q V$ per time-step. In fact, in the velocity or position version,  it is possible to merge the last substep of the
$i\rightarrow i+1$ time-step, $i = 1,\dots, I-1$ with the first substep of the subsequent time-step. This is illustrated in Table~\ref{tab:verlet} for the velocity algorithm.

There is a  feature of the velocity or position Verlet algorithms that, while not being essential for the validity of the algorithm in Table~\ref{tab:alg1} (based on preservation of volumen and reversibility), plays an important role: {\em symplecticness} \cite{arnold}, \cite{hlw}, \cite{benbook}, \cite{ssc}. When $d=1$ symplecticness is equivalent to preservation of volume (i.e.\ of planar area); when $d>1$ it is a stronger property. The symplecticness of the Verlet algorithm is a direct consequence of two facts: (i) Hamiltonian flows like $\varphi _t^A$ and  $\varphi _t^B$ are automatically symplectic and (ii) the composition of symplectic transformations is symplectic. It is well known that symplectic algorithms typically lead to energy errors smaller than its non-symplectic counterparts.

\section{Integrating the equations of motion: guidelines}
\label{sec:guide}
  The aim of this paper is to ascertain whether there exist alternative integrators that improve on the performance of the Verlet algorithm {\em within HMC and related simulations}. We limit our attention to one-step integrators where the approximation  at time $(i+1)h$ is recursively computed as $(q_{i+1},p_{i+1})= \psi_h(q_i,p_i)$. Then the transformation required by the algorithm is given
 by performing $I$ time-steps as in (\ref{eq:I}). If $\psi_h$ is volume-preserving (reversible), then $\Psi_h$ will also be volume-preserving (reversible).

The following considerations give some guidelines for the choice of integrator:
\begin{enumerate}
\item In \lq general purpose\rq\ integrations, the error after $I$ time-steps (global error)
\begin{equation}\label{eq:globalerror}
\Psi_{h,I}(q,p)-\varphi_{Ih}^H(q,p)
\end{equation}
is of paramount importance. Here we are  interested
in {\em energy errors}\footnote{While $\Delta$ depends on $h$ and $I$, this dependence is not incorporated to the notation to avoid cumbersome formulae.}
$$
\Delta(q,p) =H(\Psi_{h,I}(q,p))-H(\varphi_{Ih}^H(q,p))
$$
or, by conservation of energy,
\begin{equation}\label{eq:energyerror}\Delta(q,p) = H(\Psi_{h,I}(q,p))-H(q,p),
\end{equation} as only these determine the acceptance probability.

\item The {\em sign} of the energy error matters: $\Delta(q^{(n)},p^{(n)})<0$ always leads to acceptance of the proposal.
\end{enumerate}

In connection with the second item, it is remarkable that (see Fig.~\ref{fig:rever}) if $\Psi$ is a reversible transformation and $(q,p)$ is a point in phase space with an energy {\em increase}
$\Delta(q,p) = H(\Psi(q,p))-H(q,p) > 0$, then the point $(q^*,-p^*)$ obtained by flipping the momentum in $\Psi(q,p)$ leads
a {\em decrease} of the same magnitude
$$
\Delta(q^*,-p^*) =  H(\Psi(q^*,-p^*))-H(q^*,-p^*) = -\Delta(q,p) <0.
$$
Applying this argument to each point of a  domain $D$, we see that if the transformation is also volume-preserving, to each domain $D$ with $\Delta >0$ there corresponds a domain $S(\Psi(D))$ {\em of the same volume} with $\Delta <0$. The conclusion is that, speaking informally, for the algorithm in Table~\ref{tab:alg1} the phase space will always be divided into two regions \lq of the same volume\rq, one with $\Delta >0$ and the other with $\Delta  <0$ (and hence leading to acceptance).

It would  be wrong to infer from here that the acceptance rate should always be at least
50\%. In fact,  the standard volume (Lebesgue measure in phase space) is of little relevance and  we are rather interested  in the measure $\Pi$ in (\ref{eq:pdfpq}), as this gives the distribution of  $(q^{(n)},p^{(n)})$ at stationarity of the Markov chain. Note in Fig.~\ref{fig:rever}, that a domain
$D$ with $\Delta H>0$ as above has lower values of $H$ and carries more probability under $\Pi$ than the corresponding $S(\Psi(D))$; therefore when averaging $\Delta$ {\em with respect to $\Pi$}
the symmetry of the roles of the domains with positive and negative $\Delta$ will not be complete.

More precisely if
$$
 \mathbb{E}(\Delta)= \int_{\mathbb{R}^{2d}} \Delta(q,p) \exp\big(-H(q,p)\big)\,dq\,dp
$$
denotes the average energy error, in Fig.~\ref{fig:rever} we may observe
$$
 \mathbb{E}(\Delta) = - \int_{\mathbb{R}^{2d}} \Delta(q,p) \exp\big(-H(\Psi(q,p))\big)\,dq\,dp
$$
(an analytic proof is provided in \cite{natesh}). Thus
\begin{eqnarray*}
 \mathbb{E}(\Delta) &=& \frac{1}{2} \int_{\mathbb{R}^{2d}} \Delta(q,p) \Big[\exp\big(-H(q,p)\big)-\exp\big(-H(\Psi(q,p))\big)\Big]\,dq\,dp\\
&=&\frac{1}{2} \int_{\mathbb{R}^{2d}} \Delta(q,p) \Big[1-\exp\big(-\Delta(q,p)\big)\Big]\exp\big(-H(q,p)\big)\,dq\,dp
\end{eqnarray*}
and from here one may prove \cite{natesh}
$$
0\leq  \mathbb{E}(\Delta) \leq \int_{\mathbb{R}^{2d}} \Delta(q,p)^2 \exp\big(-H(q,p)\big)\,dq\,dp.
$$
This is a rigorous bound very relevant to our aims. It shows that the average energy error $ \mathbb{E}(\Delta)$ is of the order of $\Delta^2$ and not of the order of $\Delta$, as one may first have guessed; the result holds under the only hypotheses that the transformation $\Psi$ is volume-preserving and reversible.
If $\Psi=\Psi_{h,I}$ corresponds to an integrator of order $\nu$ , then the global error (\ref{eq:globalerror}) and the energy error (\ref{eq:energyerror}) may be bounded as $\mathcal{O}(h^\nu)$ provided that $V$ is smooth and $Ih$ remains bounded above and, accordingly, $ \mathbb{E}(\Delta) = \mathcal{O}(h^{2\nu})$ (see \cite{natesh} for technical details): {\em for our purposes the order of the method is doubled}. Reversible integrators have necessarily an even order $\nu$, Verlet has $\nu = 2$ and $ \mathbb{E}(\Delta) = \mathcal{O}(h^4)$; a fourth-order integrator would have $ \mathbb{E}(\Delta) = \mathcal{O}(h^8)$.
To sum up: due to the symmetries inbuilt in the situation, {\em average} size  of $\Delta$ will be smaller than one would have first anticipated
(see the numerical illustrations at the end of the next section).

\section{Integrating the equations of motion: the model problem}
\label{sec:model}
 A traditional approach  in the analysis of  integrators consists in the detailed  study of the application of the numerical method to the model scalar linear equation $dy/dt = \lambda y$. The conclusions  are then easily extended, via diagonalization, to general linear, constant coefficient problems and it is hoped that they  also possess  some relevance in nonlinear situations. From a negative point of view: methods that are not successful for the model equation cannot be recommended for real problems.

\subsection{The univariate case}
\label{ss:univariate}
  In our setting, a similar approach leads us to consider integrators as applied to the harmonic oscillator with Hamiltonian
  \begin{equation}\label{eq:harmonicH}
 H = \frac{1}{2}(p^2+q^2), \qquad q,p\in \mathbb{R},
 \end{equation}
  and equations of motion
\begin{equation}\label{eq:harmoscpart}
\frac{d}{dt} q = p,\qquad \frac{d}{dt} p = -q.
\end{equation}
 From the sampling point of view, this corresponds to studying the case where the target (\ref{eq:pdf}) is the standard univariate Gaussian distribution, the mass matrix is $M=1$ and (\ref{eq:pdfpq}) is a bivariate
Gaussian  with zero mean and unit covariance matrix.\footnote{We emphasize that it makes no practical sense to use a Markov chain algorithm to sample from a Gaussian distribution, just as it makes no sense to integrate numerically the equation $dy/dt = \lambda y$. In both cases it is a matter of considering simple problems as a guide to the performance of the algorithms in more realistic circumstances.} We remark that the relevance of this simple model problem to realistic quantum chromodynamics computations has been discussed in \cite{joo}.

In matrix form, the solution flow  of (\ref{eq:harmoscpart}) is given by
\begin{equation}\label{eq:rotation}
\left[ \begin{matrix}q(t)\\p(t)\end{matrix}\right] = M_t\left[ \begin{matrix}q(0)\\p(0)\end{matrix}\right],\qquad
M_t
=
\left[ \begin{matrix}\phantom{-}\cos t & \sin t\\ -\sin t & \cos t\end{matrix}\right].
\end{equation}
For all integrators of practical interest, a time-step $(q_{i+1},p_{i+1})= \psi_h(q_i,p_i)$ may be expressed as
\begin{equation}\label{eq:harmonicintegrator}
\left[ \begin{matrix}q_{i+1}\\p_{i+1}\end{matrix}\right] = \tilde{M}_h\left[ \begin{matrix}q_i\\p_i\end{matrix}\right],\qquad
\tilde{M}_h=
\left[ \begin{matrix}A_h& B_h\\ C_h & D_h\end{matrix}\right]
\end{equation}
for suitable method-dependent coefficients $A_h$, $B_h$, $C_h$, $D_h$, and the evolution over $i$ time-steps is then given by
\begin{equation}\label{eq:harmonicintegrator2}
\left[ \begin{matrix}q_{i}\\p_{i}\end{matrix}\right] = \tilde{M}_h^i\left[ \begin{matrix}q_0\\p_0\end{matrix}\right].
\end{equation}
 For a method of order $\nu$,

 \begin{equation}\label{eq:ordermatrix}
 \tilde{M}_h= M_h+\mathcal{O}(h^{\nu+1}), \qquad h\rightarrow 0,
 \end{equation}
  so that $\tilde{M}_h^i= M_{ih}+\mathcal{O}(h^\nu)$, as $h\rightarrow 0$ with $ih$  bounded above.

We restrict our interest hereafter to integrators that are both reversible and volume-preserving (symplectic since here $d=1$). For the model problem, (\ref{eq:reversible2}) leads to $A_h = D_h$ and (\ref{eq:volpre}) implies $A_hD_h-B_hC_h  = 1$. It is well known that then there are four possibilities:
\begin{enumerate}
\item $h$ is such that $|A_h| > 1$. In that case $\tilde{M}_h$ has spectral radius $>1$ and therefore the powers $\tilde{M}_h^i$ grow exponentially with $i$. For those values of $h$ the method is {\em unstable} and does not yield meaningful results.

\item  $h$ is such that $|A_h|< 1$. In that case, $\tilde{M}_h$  has complex conjugate eigenvalues of unit modulus and the powers $\tilde{M}_h^i$, $i = 0, 1,\dots$ remain bounded. The integration is then said to be {\em stable}.

\item $A_h = \pm 1$ and $|B_h|+|C_h| >0$. Then  the powers  $\tilde{M}_h^i$ grow linearly with $i$ (weak instability).

\item  $A_h = \pm 1$, $B_h=C_h=0$, i.e.\ $\tilde{M}_h = \pm I$ (stability).

\end{enumerate}

For a consistent method, $A_h = 1-h^2/2+\mathcal{O}(h^3)$, as $h\rightarrow 0$, and therefore Case 2 above holds for $h$ positive and sufficiently small. The {\em stability interval} of the method is defined as the largest interval $(0,h_{\rm max})$ such that the method is stable for each $h$, $0<h<h_{\rm max}$.

For $h$ such that $|A_h|\leq 1$, is expedient to introduce $\theta_h\in\mathbb{R}$  such that $A_h = D_h = \cos \theta_h$. For $|A_h|< 1$, we have $\sin\theta_h\neq 0$ and we may define
\begin{equation}
\label{eq:chi} \chi_h = B_h/\sin \theta_h.
\end{equation}
In terms of $\theta_h$ and $\chi_h$, the matrices in (\ref{eq:harmonicintegrator}) and (\ref{eq:harmonicintegrator2}) are then
\begin{equation}\label{eq:tildemh}
\tilde{M}_h
=
\left[ \begin{matrix}\cos \theta_h & \chi_h\sin \theta_h\\ -\chi_h^{-1}\sin \theta_h & \cos \theta_h\end{matrix}\right]
\end{equation}
and
\begin{equation}\label{eq:tildemhdos}
\tilde{M}_h^i
=
\left[ \begin{matrix}\cos (i\theta_h) & \chi_h\sin (i\theta_h)\\ -\chi_h^{-1}\sin (i\theta_h) & \cos (i\theta_h)\end{matrix}\right]
.
\end{equation}
In the (stable) case $A_h = \pm 1$, $B_h=C_h=0$, one has $\sin\theta_h=0$ and the matrix $\tilde M_h$ is of the form (\ref{eq:tildemh}) for arbitrary $\chi_h$.\footnote{Typically, it is still possible to define $\chi_{h}$ {\em uniquely} by continuity, i.e. by taking limits as $\epsilon\rightarrow 0$ in $\chi_{h+\epsilon} = B_{h+\epsilon}/\sin \theta_{h+\epsilon}$. A similar remark applies to the quantity $\rho(h)$ defined later, see (\ref{eq:ander}). Section \ref{s:choice} contains several examples of such a definition by continuity.}

From (\ref{eq:ordermatrix}) it is easily concluded that, for a method of order $\nu$, $\chi_h = 1+\mathcal{O}(h^\nu)$, $\theta_h = h + \mathcal{O}(h^{\nu+1})$ as
$h\rightarrow 0$.
By comparing the numerical $\tilde{M}_h^i$ in (\ref{eq:tildemhdos}) with the true $M_{ih}$ in (\ref{eq:rotation}), one sees that a method with $\theta_h = h$ would have no phase error: the angular frequency of the rotation of the numerical solution would coincide with the true angular rotation of the harmonic oscillator. On the other hand a method with $\chi_h = 1$ would have no energy error: the numerical solution would remain on the correct level curve of the Hamiltonian (\ref{eq:harmonicH}), i.e.\ on the circle
$p^2+q^2 = p_0^2+q_0^2$.  These considerations may  be made somewhat more precise with the help of the following well-known proposition (cf. Example 10.1 in \cite{ssc}), whose proof is a simple exercise and will not be given.

\begin{proposition} Consider a (reversible, volume-preserving) integrator (\ref{eq:harmonicintegrator}) and used with a {\em stable} value of $h$ so that
 $\tilde M_h$ may be written in the form (\ref{eq:tildemh}). Then $\psi_h = \varphi_h^{\tilde{H}_h}$,
where
$$
\tilde{H}_h = \frac{\theta_h}{2h}\left(\chi_h p^2+ \frac{1}{\chi_h}q^2\right)
$$
is the so-called modified (or shadow) Hamiltonian. In other words, one time-step of length $h$ of the numerical integrator coincides with the exact solution flow at time
$t=h$ of the Hamiltonian system with Hamiltonian function $\tilde{H}_h$.

As a consequence, $i$ time-steps of length $h$ coincide with the exact solution flow at time
$t=ih$ of the Hamiltonian system with Hamiltonian function $\tilde{H}_h$.
\end{proposition}

{\em Remark.} The existence of a modified Hamiltonian is not restricted to  harmonic problems: symplectic integrators  possess modified Hamiltonians such that the numerical solution (almost) coincides with the true solution of the modified Hamiltonian system, see e.g.\ the discussion in \cite{ssc}, Chapter 10.
\medskip

The preceding result implies that, for each fixed initial point $(q_0,p_0)$, the points $(q_i,p_i)$, $i = 1, 2,\dots$ obtained by iterating the integrator, $(q_{i+1},p_{i+1}) = \psi_h(q_i,p_i)$, lie on the level set $\tilde{H}_h(q,p) = \tilde{H}_h(q_0,p_0)$, i.e.\ on the ellipse
\begin{equation}\label{eq:ellipse}
 \chi_h p^2+ \frac{1}{\chi_h}q^2 =  \chi_h p_0^2+ \frac{1}{\chi_h}q_0^2.
\end{equation}

\begin{proposition} In the situation of the preceding proposition, for a transition over $I$ time-steps $(q_I,p_I) = \Psi_{h,I}(q_0,p_0)$, the energy error, $\Delta(q_0,p_0) = H(q_I,p_I) -H(q_0,p_0)$, may be bounded as
 $$\Delta(q_0,p_0) \leq \frac{1}{2}(\chi_h^2-1) p_0^2$$ if $\chi_h^2 \geq 1$ or as
$$\Delta(q_0,p_0) \leq \frac{1}{2}\left(\frac{1}{\chi_h^2}-1\right) q_0^2$$ if $\chi_h^2 \leq 1$.
\end{proposition}
{\em Proof.} We only deal with the first item; the other is similar. The ellipse (\ref{eq:ellipse}) has its major axis along the co-ordinate axis $p=0$ of the $(q,p)$ plane. Hence $2 H(q,p) = p^2+q^2$ attains its maximum on that ellipse if $p=0$ which implies $q^2 = q_0^2+\chi_h^2 p_0^2$. If $(q_I,p_I)$ happens to be at that maximum,
$2\Delta(q_0,p_0) = (q_0^2+\chi_h^2p_0^2) - (q_0^2+p_0^2)$.

\begin{proposition}\label{prop:average} In the situation of the preceding propositions, assume that $(q_0,p_0)$ is a random vector with distribution (\ref{eq:pdfpq}), (\ref{eq:harmonicH}).
Then the expectation of the random variable $\Delta(q_0,p_0)$ is given by
$$
\mathbb{E}(\Delta) =  \sin^2(I\theta_h)\:\rho(h),
 $$
 where
 $$
 \qquad \rho(h) = \frac{1}{2}\left(\chi_h^2+\frac{1}{\chi_h^2}-2\right)= \frac{1}{2}\left(\chi_h -\frac{1}{\chi_h}\right)^2\geq 0,
$$
and accordingly
$$
0\leq \mathbb{E}(\Delta)  \leq \rho(h).
$$
\end{proposition}
{\em Proof.} With the shorthand $c = \cos (I\theta_h)$, $s = \sin (I\theta_h)$, we may write
$$
2\Delta(q_0,p_0)  = \left(-\frac{1}{\chi_h}sq_0+cp_0\right)^2+\big(cq_0+\chi_hsp_0\big)^2-\big(p_0^2+q_0^2\big)
$$
or
$$
2\Delta(q_0,p_0)  = s^2\left(\frac{1}{\chi_h^2}-1\right)  q_0^2+ 2cs\left(\chi_h-\frac{1}{\chi_h}\right) q_0p_0
+s^2\big(\chi_h^2-1\big) p_0^2.
$$
Since $\mathbb{E}(q_0^2) = \mathbb{E}(p_0^2) = 1$ and $\mathbb{E}(q_0p_0) = 0$, the proof is ready.

\medskip

A trivial computation shows that, for $|A_h|<1$,
\begin{equation}\label{eq:ander}
\rho(h) =\frac{(B_h+C_h)^2}
{2(1-A_h^2)},
\end{equation}
a formula that will be used repeatedly in Section \ref{s:choice}.

{\em Remark.} It is relevant to note that in the last two propositions the bounds depend on $h$ but do not grow with the number $I$ of time-steps.
It is typical of symplectic integration that the energy error does not grow unboundedly as $t$ increases, see \cite{hlw}, \cite{ssc}.
\medskip

Let us illustrate the preceding results in the case of the Verlet integrator. The velocity version has $A_h=1-h^2/2$, $B_h = h$; therefore the stability interval is
$0<h<2$ (which is well known to be optimally long, see Section \ref{ss:optimal} below) and, for those values of $h$,
$$
\chi_h^2 = \frac{h^2}{1-\left(1-\frac{h^2}{2}\right)^2} = \frac{1}{1-\frac{h^2}{4}}>1.
$$
The bound in Proposition 2 reads
\begin{equation}\label{eq:boundverlet}
\Delta(q_0,p_0) \leq \frac{h^2}{8 (1-\frac{h^2}{4})}\:p_0^2.
\end{equation}
For $h = 1$, $\Delta(q_0,p_0) \leq p_0^2/6$; therefore, if $-2 < p_0 <2$ (an event that for a standard normal distribution has probability $>95\%$), then
$\Delta(q_0,p_0) < 2/3$ which results in a probability of acceptance $\geq 51\%$, regardless of the number $I$ of time-steps.

The position Verlet integrator has $\chi_h^2 = 1-h^2/4 < 1$ provided that $0<h<2$. Proposition 2 yields
$$
\Delta(q_0,p_0) \leq \frac{h^2}{8 (1-\frac{h^2}{4})}\:q_0^2
$$
(as one may have guessed from (\ref{eq:boundverlet})).

 From Proposition 3,
for both the velocity and the position versions,
\begin{equation}\label{eq:rhoverlet}
\mathbb{E}(\Delta)\leq\rho(h)= \frac{h^4}{32(1-\frac{h^2}{4})}.
\end{equation}
(Note the exponent 4 in the numerator in agreement with the discussion in the preceding section.) For $h=1$ the expected energy error is  $\leq 1/24$. Halving $h$ to $h=1/2$, leads to
an expected energy error  $\leq 1/480$.

{\em Remark.} A comparison of a given integrator (\ref{eq:harmonicintegrator}) with (\ref{eq:rotation}) shows that
\begin{equation}\label{eq:harmonicintegratorswap}
\left[ \begin{matrix}q_{i+1}\\p_{i+1}\end{matrix}\right] = \tilde{M}_h\left[ \begin{matrix}q_i\\p_i\end{matrix}\right],\qquad
\tilde{M}_h=
\left[ \begin{matrix}A_h& -C_h\\ -B_h & D_h\end{matrix}\right]
\end{equation}
is a second integrator of the same order of accuracy. The integrators (\ref{eq:harmonicintegrator}) and (\ref{eq:harmonicintegratorswap}) share the same interval of stability and the same $\theta_h$. The function $\chi_h$ of (\ref{eq:harmonicintegratorswap}) is obtained by changing the sign of the reciprocal of the function  $\chi_h$ of (\ref{eq:harmonicintegrator}). Accordingly, (\ref{eq:harmonicintegrator}) and (\ref{eq:harmonicintegratorswap}) share a common $\rho(h)$. The velocity Verlet algorithm and the position Verlet algorithm provide an example of this kind of pair of integrators.

\subsection{The multivariate case}
\label{ss:multivariate}
We now consider general Gaussian targets
$$
\pi(q)\propto \exp\left(-\frac{1}{2}q^TC^{-1}q\right)
$$
($C$ is a symmetric, positive-definite, $d\times d$ matrix of covariances). Elementary results on the simultaneous diagonalization of two quadratic form show that there is a canonical linear change of variables $q = LQ$, $p = L^{-T}Q$ that brings the Hamiltonian
$$
H = \frac{1}{2} p^TMp+\frac{1}{2}q^TC^{-1}q
$$
to the format
$$
\frac{1}{2}P^TP+Q^TDQ,
$$
where $D$ is a diagonal matrix with positive diagonal entries $\omega_j^2$. It is clear that in the new variables the equations of motion are uncoupled:
$$
\frac{d}{dt} Q_{(j)} =  P_{(j)},\qquad \frac{d}{dt} P_{(j)} =-\omega_j^2 Q_{(j)},\qquad j=1,\dots,d;
$$
in fact this uncoupling is standard in the classical theory of small oscillations around stable equilibria of mechanical systems \cite{arnold}.

The scaled variables $\bar{P}_{(j)} = P_{(j)}$, $\bar{Q}_{(j)} = \omega_j Q_{(j)}$ are uncorrelated and possess standard normal distributions. For these variables, the equations of motion  read:
\begin{equation}\label{eq:harmoscgeneral}
\frac{d}{dt} \bar{Q}_{(j)} = \omega_j \bar{P}_{(j)},\qquad \frac{d}{dt} \bar{P}_{(j)} =-\omega_j \bar{Q}_{(j)}.
\end{equation}

Now for all integrators of practical interest, the changes of variables above commute with the time-integration, i.e.\ the application of the change of variables to the numerically computed $q$, $p$ vectors yields the same results as the numerical integration of the differential system written in the new $\bar Q$, $\bar P$, variables.
 Integrating the $j$-oscillator in (\ref{eq:harmoscgeneral}) with time-step length $h$ is equivalent to integrating the standard oscillator (\ref{eq:harmoscpart}) with
 time-step length $\omega_j h$. Furthermore in the variables $\bar Q$, $\bar P$ the value of the original energy $H$ is simply
 $$
 \sum_j \frac{1}{2} \left(\bar P_{(j)}^2 +\bar Q_{(j)}^2\right).
 $$Therefore, by applying Proposition 3 to each of the individual oscillators and then summing over $j$, we conclude that, at stationarity, the error in the total energy $H$ satisfies, for stable $h$,
 \begin{equation}\label{eq:multivariate}
0\leq  \mathbb{E}(\Delta) \leq \sum_{j=1}^d \rho(\omega_j h).
 \end{equation}
 Thus the function $\rho(h)$, defined in the context of the standar harmonic oscillator is really relevant to simulations of all Gaussian measures, regardless of the choice of (symmetric, positive-definite) mass matrix.

 \subsection{Numerical illustration}
 \label{ss:illustration}

\begin{figure}
\begin{center}
\includegraphics[scale=0.50]{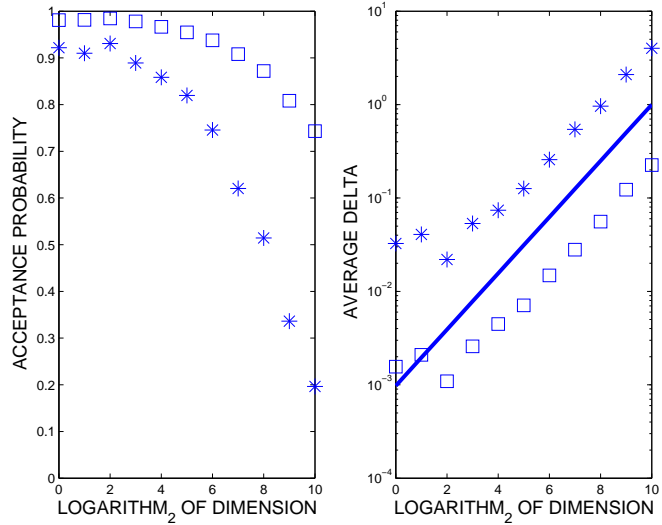}
\end{center}
\caption{\em Position Verlet algorithm. The time-step length is randomized with mean $h_0 = 1/d$ (stars) or $h_0 =1/(2d)$ (squares). On the left the observed fraction of accepted steps as a function of the number of variates $d= 1,2,\dots, 1024$. On the right the time-average of the energy increment as a function of $d$. The straight line  corresponds to
an increase proportional to $d$.}
\label{fig:halving}
\end{figure}

 We have implemented the HMC algorithm based on the position Verlet integrator for the target given by
 $$
\propto\exp\left(\frac{1}{2} \sum_{j=1}^d j^2 q_{(j)}^2\right)
 $$
 for eleven choices of the number of variates  $ d = 1,2,4,\dots , 1024$. (This distribution arises by truncating a well-known Gaussian distribution on a Hilbert space, see details in \cite{beskos}.) The mass matrix was chosen to be the identity so that the frequencies
 in the harmonic oscillators are $\omega_j = j$ and stability requires that the step-length be chosen $\leq 2/d$. The chain was started with $q^{(0)}$ at stationarity and $N=5000$ samples $q^{(n)}$ were generated. In {\em all} experiments in this paper, the step-length was randomized at the beginning of each Markov step by allowing $\pm 20\%$ variations around a mean value $h_0$
 $$ h = (1+u) h_0,\qquad u\sim \mathcal{U}(-0.2,0.2);
 $$
 among other benefits, this recipe---taken from \cite{neal}--- ensures that the observed results are not contingent on a special choice of time-step length.
 We first set $h_0 = 1/d$ (half the maximum allowed by stability). The number of time-steps was chosen as $I = 2d$ so that $T = Ih \approx 2$, a reasonable value to uncorrelate succesive samples of the \lq slowest\rq\ variate $q_{(1)}$. The results are displayed in Fig.~\ref{fig:halving}. The left panel presents the observed fraction of accepted steps; as expected (energy and energy error grow with the number of degrees of freedom) the fraction decreases as  $d$ increases and for $d= 1024$ is $\approx 20\%$. (Let us observe that, according to Table~\ref{tab:alg1}, with an energy error $\Delta(q^{(n)},p^{(n)}) = 1$ the proposal $q^*$ will be accepted with probability $\exp(-1) > 36\%$.) The figure shows that choosing $h_0$ to ensure stable integrations is not enough to achieve high rates of acceptance when the dimensionality of the problem is large. The choice $h_0 = 1/d$ works very well in this example for $d$ less than, say, 50.

 The right panel displays
 $$
 \frac{1}{N} \sum_{n=0}^{N-1}\Delta(q^{(n)},p^{(n)})
 $$
i.e.\ the observed {\em time-average}  of the
 energy error; this is seen to grow linearly with $d$ in agreement with the behavior of the {\em expectation} in (\ref{eq:multivariate}):
 $$
  \sum_{j=1}^d \rho(j h_0)= \sum_{j=1}^d \rho(j/d)\approx d\: \int_0^1 \rho(z)\,dz;
 $$
 in the language of statistical physics, the time-average and the ensemble average coincide, i.e.\ the behavior of the chain is ergodic.

Next we halved the time-step size to make it a quarter of the maximum allowed by stability ($h_0 = 1/(2d)$, $I =4d$). The Verlet integrator works well (acceptance above $70\%)$ with as many as 1,000 variates.

The right panel in the figure very clearly bears out the $h_0^4$ behavior of the average energy error. Accordingly, {\em halving the value of $h_0$ makes it possible to multiply by 16 the number of variates}. The conclusion is that, for the problem at hand and if the dimensionality is not exceptionally high, {\em the Verlet integrator may operate well even if the scaled (nondimensional) time-steps $h_0\omega_j$
are not much smaller than the upper limit imposed by stability} (say if the maximum over $j$ of $h_0\omega_j$ is between $1/2$ and $1$).

\section{Splitting methods}
In this paper we try to replace the Verlet formulae (\ref{eq:comp1}) or (\ref{eq:comp2}) by more sophisticated palindromic compositions such as:
\begin{equation}\label{eq:comp3}
\psi_{h} = \varphi_{b_1h}^B \circ \varphi_{a_1h}^A\circ\varphi_{b_2h}^B \circ \varphi_{a_2h}^A\circ\varphi_{b_2h}^B  \circ \varphi_{a_1h}^A\circ \varphi_{b_1h}^B,
\end{equation}
or
\begin{equation}\label{eq:comp4}
\psi_{h} = \varphi_{a_1h}^A \circ \varphi_{b_1h}^B\circ\varphi_{a_2h}^A \circ \varphi_{b_2h}^B\circ\varphi_{a_2h}^A  \circ \varphi_{b_1h}^B\circ \varphi_{a_1h}^A
\end{equation}
($a_j$ and $b_j$ are real parameters).
For the reasons outlined in Section 2, the single time-step mappings $\psi_{h}$  in (\ref{eq:comp3}) or (\ref{eq:comp4}) are volume-preserving, reversible and symplectic. In order to simplify the notation, we shall  use the symbols
$$
(b_1,a_1,b_2,a_2,b_2,a_1,b_1)
$$
and
$$
(a_1,b_1,a_2,b_2,a_2,b_1,a_1)
$$
to refer to (\ref{eq:comp3}) and (\ref{eq:comp4}) respectively.

With a similar notation, one may consider $r$-stage compositions, $r = 1,2,\dots$:
\begin{equation}\label{eq:comp5}
\stackrel{2r+1  \ {\rm  letters}}{\overbrace{(b_1,a_1,b_2,\dots, a_1,b_1)}}
\end{equation}
or
\begin{equation}\label{eq:comp6}
\stackrel{2r+1  \ {\rm  letters}}{\overbrace{(a_1,b_1, a_2,\dots, b_1,a_1)}}.
\end{equation}
Obviously (\ref{eq:comp6}) requires $r$ evaluations of $\nabla_q V$ at each time-step. The same is essentially true for (\ref{eq:comp5}), because, as discussed for the velocity Verlet algorithm, the last evaluation of $\nabla_q V$ at the current time-step is re-used at the next time-step. As in the Verlet algorithms, both (\ref{eq:comp5}) and (\ref{eq:comp6}) are best implemented by combining the last substep of the current time-step with the first substep of the subsequent time-step.

\subsection{Taylor expansion of the energy error}
The Lie bracket of vector fields plays an important role in the analysis of splitting integrators (see \cite{acta} and \cite{03}). In the Hamiltonian context, the vector-valued Lie bracket may be advantageously replaced by the real-valued Poisson bracket of the Hamiltonian functions; recall that if $F$ and $G$ are smooth real-valued functions in phase space, their Poisson bracket is, by definition, the function \cite{arnold}, \cite{ssc}
$$
\{F,G\} = \sum_{j=1}^d \left(\frac{\partial F}{\partial q_{(j)}}\frac{\partial G}{\partial p_{(j)}}-
\frac{\partial F}{\partial p_{(j)}}\frac{\partial G}{\partial q_{(j)}}\right)
$$
(as before, $q_{(j)}$ and $p_{(j)}$ are the scalar components of the vectors $q$ and $p$).
The properties of (\ref{eq:comp5}) or (\ref{eq:comp6}) are encapsulated in the corresponding modified Hamiltonian, which, for consistent methods and in the limit $h \rightarrow 0$, has an expansion
\begin{eqnarray}
  \tilde H_h & = & H + h^2 k_{3,1} \{A,A,B\} + h^2 k_{3,2} \{B,A,B\} \nonumber   \\
  & & +h^4 k_{5,1} \{A,A,A,A,B\}   +h^4 k_{5,2} \{B,A,A,A,B\}\nonumber\\&& + h^4 k_{5,3} \{A,A,B,B,A\} + h^4 k_{5,4} \{B,B,A,A,B\} + \mathcal{O}(h^6),
  \label{eq:tildeH}
\end{eqnarray}
where $k_{\ell,m}$ are polynomials in the coefficients $a_j$, $b_j$ and expressions like $\{A,A,B\}$ (or  $\{A,A,A,A,B\}$) are abbreviations to refer to iterated Poisson brackets
$\{A,\{A,B\}\}$ (or  $\{A,\{A,\{A,\{A,B\}\}\}\}$).\footnote{Here we have used that  $\{B,B,A,B\}=0$, a condition that is implied by the fact that $A$ is quadratic in the momentum
$p$ (see (\ref{eq:kpot})). For splittings with $\{B,B,A,B\}\neq 0$ there are six $\mathcal{O}(h^4)$ terms in the expansion of $\tilde H_h$.} Order $\nu\geq 4$ is then equivalent to the conditions
$ k_{3,1}= k_{3,2}=0$, while order $\nu\geq 6$ would require, in addition, $ k_{5,1}= k_{5,2}=k_{5,3}=k_{5,4}=0$.

By using the Lie formalism,  the Taylor expansion of the energy after one time-step is found to be (\cite{ssc}, Section 12.2)
\begin{eqnarray*}
H(q_{i+1},p_{i+1}) &=& \exp(-h{\cal L}_{\tilde H_h})H(q_i,p_i)\\& = &H(q_i,p_i) - h {\cal L}_{\tilde H_h}H(q_i,p_i)
+ \frac{1}{2}h^2 {\cal L}_{\tilde H_h}^2H(q_i,p_i)+\cdots,
\end{eqnarray*}
where ${\cal L}_{\tilde H_h}$ is the Lie operator
$
{\cal L}_{\tilde H_h}(\cdot) = \{\tilde H_h,\cdot\}
$. A trite computation then yields
\begin{equation}\label{eq:expansiondelta}
\Delta(q_i,p_i) =  h^3 k_{3,1} \{A,A,A,B\} + h^3 (k_{3,1} + k_{3,2}) \{A,B,A,B\}  + \mathcal{O}(h^4)
\end{equation}
(the iterated brackets in the right-hand side are evaluated at $(q_i,p_i)$). Thus, when $h$ is small, $E^*=k_{3,1}^2+(k_{3,1} + k_{3,2})^2$ is a measure of energy errors.
The velocity Verlet integrator has a value of $E^*$ larger than that of its position counterpart.

\subsection{Optimal stability of the Verlet integrator}
\label{ss:optimal}

The application  of a method of the form (\ref{eq:comp5}) or (\ref{eq:comp6}) to the standard harmonic oscillator (\ref{eq:harmonicH}) results in a recursion of the form
(\ref{eq:harmonicintegrator}) (of course $A_h=D_h$, $A_h^2-B_hC_h=1$ due to reversibility and volume preservation). Additionally $A_h$ is a  polynomial of degree $\leq r$ in $\zeta = h^2$ and, for consistent methods, $A_h = 1-h^2/2+O(h^4)$ as $h\rightarrow 0$. By using well-known properties of the Chebyshev polynomials it is not difficult to prove that a polynomial $P(\zeta)$ of degree $\leq r$ subject to the  requirements $P(0) = 1$, $P^\prime(0)=-1/2$ cannot satisfy $-1\leq P(\zeta)\leq 1$ for $0<\zeta< \zeta_{\rm max}$ if $\zeta_{\rm max} > 4 r^2$.
 This proves that there is no choice of coefficients for which the stability interval $(0,h_{\rm max})$ of (\ref{eq:comp5}) or (\ref{eq:comp6}) has $h_{\rm max}> 2r$ (see \cite{jeltsch}).\footnote{
 By arguing as in \cite{spijker}, the conclusion $h_{\rm max}\leq 2r$ also follows from the well-known Courant-Friedrichs-Lewy restriction for the integration of hyperbolic partial differential equations. Consider the familiar wave equation $\partial_t Q(x,t) = P(x,t)$, $\partial_t P(x,t) = \partial_x^2Q(x,t)$ with periodic boundary conditions and discretize the space variable $x$ by standard central differences. The highest frequency is $\omega = 2/\Delta x$. A consistent, explicit, one-step integrator using $r$ force evaluations per time-step with stability interval longer than $2r$ would yield a convergent approximation to the wave problem for $h/\Delta x > r$ and this violates the CFL restriction.
 } 
 Furthermore
 since the velocity  Verlet algorithm has stability interval $0<h<2$, the  concatenation $\psi_h = \psi^{\rm VelVer}_{h/r}\circ \cdots \circ \psi^{\rm  VelVer}_{h/r}$ of $r$ time-steps of length $h/r$ is a method of the form
 (\ref{eq:comp5}) that attains the optimal value $h_{\rm max} = 2r$; similarly the $r$-fold concatenation $\psi_h = \psi^{\rm  PosVer}_{h/r}\circ \cdots \circ \psi^{\rm  PosVer}_{h/r}$ is a method of the form (\ref{eq:comp6}) with optimal stability interval.

 When comparing the size of stability intervals the computational effort has to be taken into account: with a given amount of computational work, an integrator with fewer function evaluations per time-step may take shorter time-steps to span a given time interval $0\leq t\leq T$. It is therefore a standard practice to {\em normalize} the length $h_{\rm max}$ of the stability interval of explicit integrators by dividing by the number of force evaluations per time-step. According to the preceding discussion, the (position or velocity) Verlet algorithm and its concatenations have an optimal normalized stability interval of length 2.
Integrators with short normalized stability intervals  are of no interest here as they cannot compete with the Verlet scheme (see the conclusion at the end of the preceding section). In particular and as we shall see later, high-order methods proposed in the literature have stability intervals far too short and cannot compete in practice with the performance of the Verlet scheme in HMC simulations.

\section{Choice of coefficients}
\label{s:choice}
In this section we address the question of how best to choose the number of stages $r$ and the coefficients $a_j$ and $b_j$ in (\ref{eq:comp5}) or (\ref{eq:comp6}). In the derivation of numerical integrators, both for  general and HMC use (see e.g. \cite{takaishi}), it is customary to first determine $r$ to achieve a target order of accuracy $\nu$ and  to then use any remaining free parameters to minimize the  error constants. In the Hamiltonian scenario a standard way of minimizing the error constants is to reduce the coefficients of the modified Hamiltonian (\ref{eq:tildeH}). For instance, for a method of second order, one would try to minimize  some norm of the vector $(k_{3,1},k_{3,2})$. However  the ideas of order of accuracy and local error constants  both refer to the asymptotic behavior of the integrator as $h\rightarrow 0$ and we have seen in Sections \ref{sec:guide} and \ref{sec:model} that the Verlet integrator is capable of performing well in HMC simulations for rather large values of the time-step $h$. Accordingly, we shall determine $r$, $a_j$ and $b_j$ by means of a different strategy based on Gaussian models. Given a family of methods, we shall express the quantity $\rho(h)$ defined in Proposition~\ref{prop:average} as a function of the method coefficients and then we shall choose these coefficients to minimize
$$
\| \rho \|_{(\bar h)} = \max_{0<h<\bar h} \rho(h),
$$
where $\bar h$ is a suitable maximum time-step. (It is tacitly understood that $\bar h$ is {\em smaller} than the maximum step-size allowed by stability.) More precisely, since
in Section \ref{sec:model} we saw that, for the standard harmonic oscillator, the Verlet method  is capable of performing well in HMC simulations when $h\approx 1$, an efficient $r$-stage method should be able to operate well  with $h\approx r$ (if that were not the case, use of the Verlet method with time-step $h/r$ would outperform the more complex integrator with time-step $h$). Following this rationale,
we set $\bar{h} = r$ and use $\| \rho \|_{(r)}$  as a metric for the quality of an integrator within the HMC algorithm.\footnote{It is clear that it is also possible to consider alternative values of $\bar h$ or norms different from the maximum norm. Such a fine tuning will not be undertaken here.}

\subsection{Two-stage methods: $
(a_1,b_1,a_2,b_1,a_1)
$}

We start by discussing in detail methods of the form
$
(a_1,b_1,a_2,b_1,a_1)
$. Consistent integrators have to satisfy $b_1 = 1/2$, $a_2 = 1-2a_1$ and this leaves the one-parameter family
 \begin{equation}\label{eq:family}
 (a_1,1/2, 1-2a_1,1/2,a_1).
 \end{equation}
 The choices $a_1=0$ and $a_1 = 1/2$ are singular; for them the integrator reduces to the velocity Verlet and position Verlet algorithm respectively. Furthermore for  $a_1= 1/4$ one time-step   $\psi_h$ of (\ref{eq:family}) coincides with the concatenation $\psi^{\rm PosVer}_{h/2}\circ \psi^{\rm PosVer}_{h/2}$ of two time-steps of length $h/2$ of the position Verlet integrator. A standard computation (say using the Baker-Campbell-Hausdorff\footnote{The use of the BCH formula to analyze splitting algorithms may be bypassed by following the approach in \cite{ass}.}  formula \cite{ssc}) yields
 $$
k_{3,1} = \frac{12a_1^2-12a_1+2}{24},\qquad k_{3,2} = \frac{-6a_1+1}{24}.
$$
Since  no choice of $a_1$  leads to $ k_{3,1}= k_{3,2}=0$, no method of the family achieves order 4, see (\ref{eq:tildeH}). The expression $E=k_{3,1}^2+k_{3,2}^2$ that measures the leading error terms  turns out to be a convex function of the free parameter $a_1$ and has a minimum value $E\approx 7\times 10^{-5}$ at $a_1 \approx 0.1932$, as first observed by McLachlan \cite{mclachlan}.

For comparison, $a_1 = 1/4$, which is equivalent to the standard position Verlet integrator, yields a much worse $E \approx 9 \times 10^{-3}$. Therefore the choice $a_1 \approx 0.193$ has been recommended in the HMC context\footnote{The paper \cite{takaishi}  does not cite McLachlan \cite{mclachlan} and attributes the method to later papers by Omelyan and his coworkers.} \cite{takaishi} (this paper and \cite{mclachlan} give a representation of $a_1$ in terms of surds).\footnote{The expansion (\ref{eq:expansiondelta}) may suggest to minimize $E^* = k_{3,1}^2+(k_{3,1}+k_{3,2})^2$. This leads to
$a_1 = 0.1956\dots$. We shall not be concerned with this value of $a_1$, as the method is very similar to the one derived via minimization of $E$.}

As discussed above, we here follow a different strategy, based on Gaussian models. We first find from (\ref{eq:ander})
\begin{equation}\label{eq:rhodos}
  \rho(h) = \frac{h^4 \big(2a_1^2(1/2-a_1)h^2+4a_1^2-6
   a_1+1\big)^2}
   {8 \big(2-a_1 h^2\big)
   \big(2-(1/2-a_1) h^2\big)
   \big(1-a_1(1/2-a_1) h^2\big)}.
\end{equation}
Stability is equivalent to the positivity of the denominator. Note that for $a_1=0$ or $a_1 =1/2$ the quotient (\ref{eq:rhodos}) reduces to (\ref{eq:rhoverlet}), as it should. When $a_1\leq 0$ or $a_1 \geq 1/2$, $k_{3,1}$ and $k_{3,2}$ are too large and the stability interval too small. Therefore,
 useful methods have $0< a_1 < 1/2$. In this parameter range, the stability interval is
 \begin{equation}\label{eq:stabfamily}
 0<h< \min\Big\{\sqrt{2/a_1}, \sqrt{2/(1/2-a_1)}\Big\},
 \end{equation}
 {\em provided} that $a_1 \neq 1/4$. When $a_1 = 1/4$, the product $\big(2-a_1 h^2\big)
   \big(2-(1/2-a_1) h^2\big)$ in the denominator of (\ref{eq:rhodos}) is a factor of the numerator and a simplification takes place: the fraction reduces to (\ref{eq:rhoverlet}) with $h$ replaced by $h/2$ and the stability interval is $0<h<4$ in lieu of the shorter interval $0<h<2\sqrt{2}$ in (\ref{eq:stabfamily}). This corresponds to the earlier observation that for $a_1=1/4$, the method coincides with the concatenation
   $\psi^{\rm PosVer}_{h/2}\circ \psi^{\rm PosVer}_{h/2}$.

\begin{figure}
\begin{center}
\includegraphics[scale=0.50]{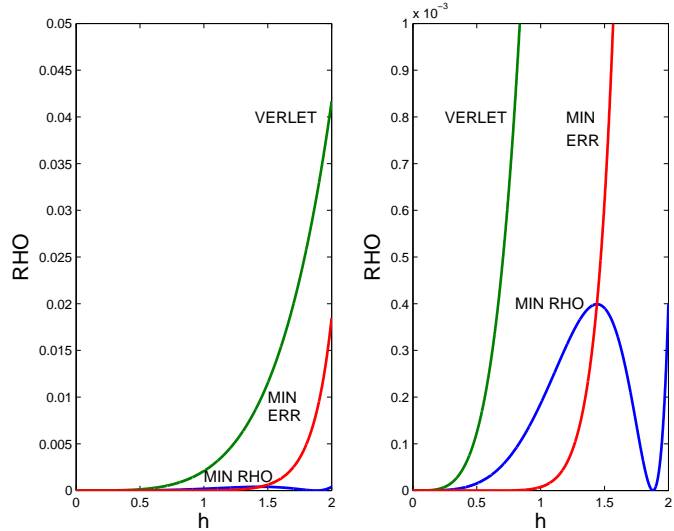}
\end{center}
\caption{\em $\rho(h)$ as a function of $h$, $0<h<2$, for three members of the family (\ref{eq:family}): $a_1= 0.25$ (equivalent to the position Verlet algorithm), $a_1$ with mimimum error coefficient $E$ ($a_1 \approx 0.1932$) and $a_1$ in (\ref{eq:nuestro}). With $a_1= 0.25$, $\rho(h)\uparrow \infty$ as $h\rightarrow 4$; for the other two choices the vertical asymptote is located at $h\approx 2.55$ and $h\approx 2.63$ respectively. The right pannel, with an enlarged vertical scale, $0< \rho< 0.001$, shows the superiority of the minimum error coefficient method for small $h$.}
\label{fig:rho}
\end{figure}

   The next task is to determine $a_1$ to minimize $\|\rho\|_{(2)}$. This yields $a_1 = 0.21178\dots$ but, to avoid cumbersome decimal expressions, we shall instead use the approximate value\footnote{For this choice of $a_1$, $k_{3,1}=0$,  so that the literature \cite{mclachlansmall} has suggested this value for cases where in (\ref{eq:tildeH}) $\{A,A,B\}$ is much larger than $\{B,A,B\}$.}
   \begin{equation}\label{eq:nuestro}
   a_1 = \frac{3-\sqrt{3}}{6} \approx 0.21132
   \end{equation}
   which gives $\|\rho\|_{(2)} \approx 5\times 10^{-4}$. For comparison, $a_1=1/4$ has a substantially larger  $\|\rho\|_{(2)} \approx 4\times 10^{-2}$ and the method of McLachlan with minimum error constant has $\|\rho\|_{(2)} \approx 2\times 10^{-2}$. Thus, when using $\|\rho\|_{(2)}$ as a metric, the minimum error-constant method provides only a marginal improvement on Verlet. In Fig.~\ref{fig:rho},  we see that, while the minimum error constant method leads to the smallest values of $\rho(h)$ for $h<1$,
   the choice (\ref{eq:nuestro}) ensures a much better behavior over the target interval $0<h<2$.

\begin{figure}
\begin{center}
\includegraphics[scale=0.40]{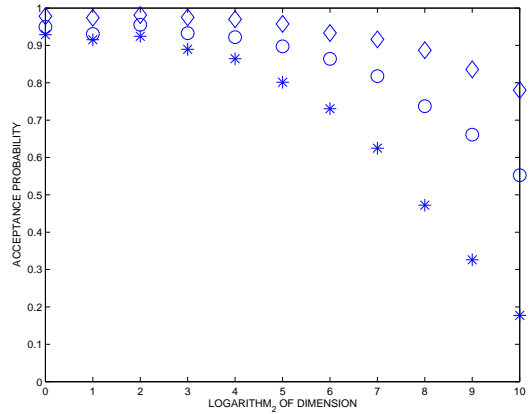}
\end{center}
\caption{\em  Observed fraction of accepted steps as a function of the number of variates $d=1,2,\dots,1024$. Stars: position Verlet with (average) step-size  $h_0 = 1/d$. Circles: two-stage method with minimum error constant, average step-size $h_0 = 2/d$. Diamonds: suggested two-stage method (\ref{eq:nuestro}), (average) step-size $h_0 = 2/d$.}
\label{fig:dosetapas}
\end{figure}

We have considered again the experiment in Section~\ref{ss:illustration}, this time comparing the position Verlet algorithm with $h_0 = 1/d$ and $I=2d$ (i.e.\ the parameters for the run marked by stars in Fig.~\ref{fig:halving}) with members of the family (\ref{eq:family}) with $h_0 = 2/d$, $I=d$, so as to equalize work. The results are shown in Fig.~\ref{fig:dosetapas}. In this problem the minimum error-constant method provides an improvement on Verlet, but its performance is markedly worse than that of the method with the value (\ref{eq:nuestro}) suggested here.  While the advantage of the method (\ref{eq:nuestro}) over Verlet in Fig.~\ref{fig:dosetapas} occurs for all values of $d$, it becomes more prominent as $d$ increases.

A comparison of Figs.~\ref{fig:halving} and \ref{fig:dosetapas} shows that the fraction of accepted steps is larger for the two-stage method (\ref{eq:nuestro}) with $h_0 = 2/d$
(diamonds) than for the Verlet algorithm with $h_0 = 1/(2d)$ (squares) (the latter simulation is twice as costly).

\subsection{Two-stage methods: $
(b_1,a_1,b_2,a_1,b_1)
$}
Let us now turn the attention to the format
$
(b_1,a_1,b_2,a_1,b_1)
$. By the remark at the end of Section~\ref{ss:univariate} the function $\rho(h)$ of such a method coincides with that of the method of the format $
(a_1,b_1,a_2,b_1,a_1)
$
based on the same sequence of numerical values of the coefficients. This leads to the integrator
$$
   b_1 = \frac{3-\sqrt{3}}{6},\quad a_1 = a_2 = 1/2, \quad b_2 = 1/2-b_1,
$$
that has
$$
k_{3,1} = \frac{6\lambda-1}{24},\qquad k_{3,1}+k_{3,2} = \frac{-12\lambda^2+18\lambda-3}{24},\qquad \lambda =  \frac{3-\sqrt{3}}{6}.
$$
For (\ref{eq:family}) with $a_1$ given in (\ref{eq:nuestro})\footnote{The lack of symmetry of the $(a_1,b_1,a_2,b_1,a_1)$ and $(b_1,a_1,b_2,a_1,b_1)$ formats is due to the fact that (\ref{eq:expansiondelta}) is not symmetric in $A$ and $B$, which in turn is a consequence of the fact that $\{B,B,A,B\}= 0$ while nothing can be said in general
about $\{A,A,B,A\}$. }
$$
k_{3,1} = 0,\qquad k_{3,1}+k_{3,2} = \frac{-12\lambda^2+18\lambda-3}{24};
$$
we find no reason to prefer in this context the $(b_1,a_1,b_2,a_1,b_1)$
sequence,
since both methods share the value of $k_{3,1}+k_{3,2}$ and $|k_{3,1}|$ is smaller for the $(a_1,b_1,a_2,b_1,a_1)$ format.

For reasons of brevity  we shall not consider again in what follows formats beginning with the letter $b$.

\subsection{Three stages}
For the three-stage format $
(a_1,b_1,a_2,b_2,a_2,b_1,a_1)
$, consistency requires $a_2 = 1/2-a_1$, $b_2 = 1-2b_1$ and therefore we have to consider the two-parameter family
$$
(a_1, b_1, 1/2-a_1, 1-2b_1, 1/2-a_1, b_1,a_1).
$$
The choice $a_1 =(1/2) (2-2^{1/3})^{-1}$, $b_1 =2a_1$ leads to a fourth-order integrator ($k_{3,1} = k_{3,2} = 0$) that goes back to Suzuki, Yoshida and others, see e.g.\ \cite{ssc}, Chapter 13.
Its stability interval is {\em very} short: approximately $0<h<1.573$.  In the HMC context, this fourth order integrator has been considered by the physics literature, see e.g. \cite{joo}, \cite{takaishi}; the former reference notes the poor stability properties.

According to our methodology, we choose the free parameters so as to minimize the maximum of $\rho(h)$ over the interval
$0<h< 3$. The situation is somewhat delicate, as we shall explain presently. Let us first consider the choice $a_1=1/6$, $b_1 = 1/3$, leading to the concatenation $\psi^{\rm PosVer}_{h/3}\circ \psi^{\rm PosVer}_{h/3} \circ \psi^{\rm PosVer}_{h/3}$, that as discussed in Section \ref{ss:optimal} possesses optimal stability interval $(0,6)$. At $h=3$, this method has  $A_h = -1$, $B_h = C_h = 0$.
 Furthermore $h=3$ is a simple root of the equations $B_h=0$ and $C_h=0$ and a {\em double} root of the equation $A_h=-1$. By the implicit function theorem,
 when the coefficients $a_1$ and $b_1$ are perturbed away from  $a_1=1/6$, $b_1 = 1/3$, the root $h=3$ of the equation $B_h=0$  moves to a  location $h_B(a_1,b_1)\approx 3$. In a similar manner, the root of $C_h=0$ moves to a location $h_C(a_1,b_1)\approx 3$, that, generically, does {\em not} coincide with $h_B(a_1,b_1)$. Now, the relation $A^2_h=1+B_hC_h$, that follows from conservation of volume, ensures that both
$h_B$ and $h_C$ are roots of $A_h = -1$. In other words,  perturbations generically change the double root $h=3$ of $A_h = -1$   present in the concatenated Verlet method into two {\em real} simple roots $h_B$, $h_C$; in the neighborhood of  such simple roots $A_h$ cannot remain $\geq -1$ and we conclude that, for  generic perturbations, the integrator is unstable near $h\approx 3$.

 In order to identify integrators (not necessarily close to the concatenated Verlet method) that do not turn unstable for $h\approx 3$ due to $A_h$ becoming $<-1$ we proceed as follows. We write $A_h$, $B_h$, $C_h$ in terms of the parameters $a_1$, $b_1$ (the expressions are cumbersome and will not be reproduced here), fix a value $\widehat h$ and consider the system of two (nonlinear) equations
\begin{equation}\label{eq:system}
A_{\widehat h} = -1,\qquad B_{\widehat h}+C_{\widehat h}=0,
\end{equation}
for the two unknowns $a_1$, $b_1$. When these relations hold, from $A_h^2-B_hC_h= 1$, we infer that $B_{\widehat h} = C_{\widehat h}=0$ and $\widehat h$ is a stable value. Furthermore, if ${}^\prime$ denotes differentiation with respect to $h$,
$$
2A_h A_h^\prime- B_h^\prime C_h-B_hC_h^\prime = 0
$$
and therefore $A_{\widehat h}^\prime=0$, so that $A_h$ will have a minimum at $h=\widehat h$ and thus remain $\geq -1$ in the neighborhood of $\widehat h$.
Note also that in (\ref{eq:ander}), the zero of the denominator at $h=\widehat h$ may be simplified with the corresponding zero of the numerator (an occurrence we already found
when discussing (\ref{eq:rhodos}) and entails an \lq enlargement\rq\ of the stability interval).

We solve the system of equations (\ref{eq:system}) and find the following family of integrators parameterized by the location $\widehat h$ of the double root of $A_h=-1$:
\begin{equation}\label{eq:oneparameter}
a_1 =\frac{1}{2} -\frac{3}{\widehat{h}^2}\pm \frac{\sqrt{9-\widehat{h}^2}}{\widehat{h}^2},\qquad
b_1 = \frac{3}{\widehat{h}^2}\pm \frac{\sqrt{9-\widehat{h}^2}}{\widehat{h}^2},\qquad 0< \widehat h\leq 3.
\end{equation}
For $\widehat h=3$,  the integrator is the concatenation of three Verlet substeps discussed above.  For $\widehat h = 2\sqrt{2}$ and the positive value of the square root,
we find $a_1=0$, $b_1=1/4$ and the integrator is $\psi^{\rm VelVer}_{h/2}\circ \psi^{\rm VelVer}_{h/2}$, with stability interval $0<h<4$. The negative value of the square root
leads to $a_1=1/4$, $b_1=1/2$ ($b_2=0$), i.e.\ to $\psi^{\rm PosVer}_{h/2}\circ \psi^{\rm PosVer}_{h/2}$,  whose stability interval is again $0<h<4$.\footnote{Not all members of the family of methods (\ref{eq:oneparameter}) are stable for all values of $h$, $0<h<3$, as instability may also occur by $A_h$ becoming larger than 1.}

Finally we determine  $\widehat h$ in (\ref{eq:oneparameter}) by minimizing $\| \rho\|_{(3)}$. This yields the parameter values
\begin{equation}\label{eq:minrhotres}
a_1=  0.11888010966548,\qquad
b_1 = 0.29619504261126,
\end{equation}
with $\| \rho\|_{(3)}= 7\times 10^{-5}$ and stability interval of length $\approx 4.67$ (the double root is located at $\widehat h \approx  2.98$).\footnote{The perturbation argument presented above for the concatenated Verlet method applies to perturbations of any member of the family (\ref{eq:oneparameter}): generic perturbations turn the double root into two simple real roots, which leads to instability near $\widehat h$.}

\begin{figure}
\begin{center}
\includegraphics[scale=0.40]{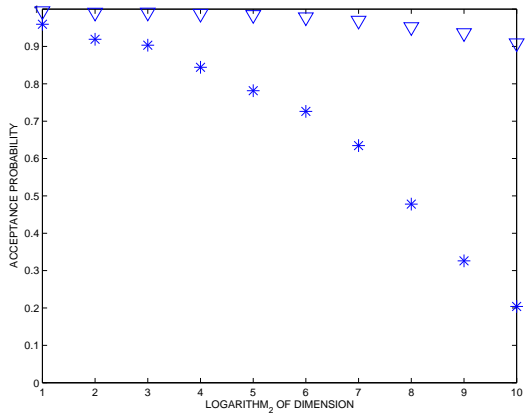}
\end{center}
\caption{\em  Observed fraction of accepted steps as a function of the number of variates $d=2,4,\dots,1024$. Stars: position Verlet algorithm with (average) step-size $h_0 = 1/d$. Triangles: three-stage method (\ref{eq:minrhotres}) with $h_0= 3/d$.}
\label{fig:tres}
\end{figure}

An illustration of the performance of the new integrator may be seen in Fig.~\ref{fig:tres} which refers again to the experiment in Section~\ref{ss:illustration}. The position Verlet algorithm is run with $h_0=1/d$ (as in Figs.~\ref{fig:halving} and  \ref{fig:dosetapas}) and, in order to equalize work, the three-stage method (\ref{eq:minrhotres}) was used with $h_0=3/d$. The number of time-steps was taken to be (the integer closest to) $2d/3$ for (\ref{eq:minrhotres}) and thrice that number for Verlet, so that the (average) final time is $T=2$ and both methods use the same number of force evaluations. The advantage of the three-stage method over both Verlet and the optimized  2-stage integrator is clearly felt.

It is of interest to point out that with the present  choice  $h_0=3/d$ the fourth-order, three-stage integrator is unstable; in fact $h_0$ would have to be halved to barely ensure stability. However, when the step-length is halved, Verlet  delivers satisfactory acceptance rates as we saw in Fig.~\ref{fig:halving}. We conclude  that the  benefits of high order only take place when $h_0$ is too small for the goals of the integration.

\subsection{Four stages}

With four stages $(a_1,b_1,a_2,b_2,a_3,b_2,a_2,b_1,a_1)$, by using a similar procedure we find the method, with $\|\rho\|_{(4)} \approx 7\times 10^{-7}$,
\begin{eqnarray}
a_1 &=& 0.071353913450279725904,\nonumber \\
a_2 &=& 0.268548791161230105820,\nonumber\\
b_1 &=& 0.191667800000000000000 \label{eq:minrhocuatro}
\end{eqnarray}
(the remaining parameter values are determined by consistency, i.e.\ $b_2 = 1/2-b_1$, $a_3 = 1-2a_1-2a_2$). The method has a stability interval of length $\approx 5.35$ (the equation $A_h = -1$ possesses a double root at $\approx 3.04$). For the target in Section~\ref{ss:illustration} with $h_0 = 4/d$ and
$I = d/2$ (which involves the same computational effort as the Verlet runs marked by stars in Figs.~\ref{fig:halving} and \ref{fig:dosetapas}) the observed fraction of accepted steps remains above $98\%$ for all values of $d= 2,2^2,\dots, 2^{10}$. Such large acceptance rates would be most welcome in variants of HMC, including the generalized HMC of Horowitz \cite{horo}, \cite{elena}, where rejections are particularly troublesome.

\section{A small molecule}

  A detailed benchmarking of the various integrators in different application examples will be considered elsewhere and is not within the scope of our work here. However, since our methodology is based on a Gaussian model problem, it is of clear interest to run some proof-of-concept experiments with non-Gaussian targets.
We have used as a test problem the Boltzmann distribution of a pentane molecule, as in the numerical comparisons in \cite{cances}. The model has fifteen degrees of freedom (the cartesian coordinates of the five carbon atoms); it includes very strong forces associated with the  carbon-carbon covalent bond length, softer forces associated with the bond and dihedral angles and also Van der Waals interactions. The number of vibrational degrees of freedom, nine as there are six corresponding to rigid-body motions, is modest and therefore we may expect that the Verlet algorithm may be able to work with step-sizes not much smaller than the maximum allowed by stability; thus the choice of problem may be considered  to be  biased in favor of Verlet.
The molecule has several stable configurations (minima of the potential energy) and therefore the target distribution is multimodal; the highly nonlinear Hamiltonian dynamics moves the molecule among the different configurational energy basins. Some degrees of freedom (i.e.\ bond lengths) have very small variances, other (such as dihedral angles) vary by substantial amounts. We set the molecule parameters as in \cite{cances}
and the results reported here correspond to an inverse temperature $\beta = 1/2$. The simulation starts from the most stable configuration and from there takes 200 Markov burn-in HMC steps to bring the chain to stationarity; after that, samples are taken from the next 512 Markov steps. Five integrations are considered:
\begin{itemize}
\item Position Verlet integrator, $h_0= 0.02$, $I = 24$.
\item Minimum error-constant two-stage integrator, $h = 0.04$ , $I =12$.
\item Two-stage  integrator (\ref{eq:nuestro}), $h = 0.04$ , $I =12$.
\item Three-stage integrator (\ref{eq:minrhotres}), $h= 0.06$, $I=8$.
\item Four-stage integrator (\ref{eq:minrhocuatro}), $h=0.08$, $I=6$.
\end{itemize}
The values $h_0= 0.02$, $I = 24$  were tuned to provide a good performance of the Verlet algorithm (performance was measured by the efficiency in computing the probabilities that the molecule is in its different configurational basin). After that, the values of $h_0$ and $I$ for the other integrators were determined to ensure that all integrations share a common computational effort.

\begin{table}
\begin{center}
\begin{tabular}{|l|r|r|}
Integrator
& $\mu$ & $\sigma$\\\hline
(One-stage) Verlet & $85\%$ & $2.0\%$\\
Two-stage minimum error-constant &$80\%$ & $1.9\%$\\
Two-stage minimum $\rho$ (\ref{eq:nuestro}) & $92\%$ & $1.2\%$\\
Three-stage minimum $\rho$ (\ref{eq:minrhotres})& $97\%$& $0.7\%$\\
Four-stage minimum $\rho$ (\ref{eq:minrhocuatro})&$97\%$ &$0.8\%$\\
\hline
\end{tabular}
\end{center}
\caption{\em Pentane molecule. Mean value and standard deviation, over 100 realizations of the Markov chain, of the observed acceptance ratio.}
\label{t:molecule}
\end{table}

 For each integrator we computed 100 realizations of the Markov chain;  Table~\ref{t:molecule} displays the mean value (over the 100 samples) of the empirical acceptance rate (after burn-in) and the associated standard deviation.
It is apparent in Table~\ref{t:molecule}
that the performance of the minimum error constant integrator is worse than that of the Verlet algorithm. No doubt this is due to the fact that the time step-sizes involved are too large for the Taylor expansions (\ref{eq:tildeH}) and (\ref{eq:expansiondelta}) to be meaningful for the problem under consideration. In fact, additional experiments with the pentane molecule prove that when the integrator with minimum error constant and Verlet are used with step-lengths that equalize work, the energy errors of the former improve on those of Verlet only for step-sizes so small that the acceptance rate for Verlet is very approximately 100\%.
On the other hand, Table~\ref{t:molecule} reveals that the two-stage integrator suggested here does
improve on the Verlet integrator. The most efficient integrations are afforded by the three and four stage schemes, even though, as pointed out before, the low-dimensionality of the problem biases this model problem against the more sophisticated integrators.

\section{Conclusions}
We have suggested a methodology for constructing efficient methods for the numerical integration of the Hamiltonian differential equations that arise in HMC and related algorithms. The new approach is based on optimizing the behavior of a function $\rho(h)$ over a relevant range of values of the step-length $h$. We have constructed new split-step integrators with two, three or four function evaluations  per time-step. Unlike integrators derived by minimizing the size of error constants, the splitting formulae suggested here are more efficient than the standard Verlet method, specially if the number of dimensions is high.

The detailed benchmarking of the new integrators will be the subject of subsequent work.

\bigskip

{\bf Acknowledgments.} The authors are thankful to A. Murua for some useful discussions. They have been supported by projects MTM2010-18246-C03-01 (Sanz-Serna)
and MTM2010-18246-C03-02  (Blanes and Casas) from Mi\-nisterio de Ciencia e
Innovaci\'on, Spain. Additionally, Casas has  been partly supported by grant NPRP 5-674-1-114 from the Qatar National Research Fund.
\end{document}